\documentclass[twocolumn,10pt]{IEEEtran}
\usepackage[bottom=0.68in,top=0.6in,left=0.625in,right=0.625in]{geometry}
\usepackage{ifpdf, flushend,subfigure,tabularx}
\usepackage{graphicx}
\usepackage{epstopdf}
\graphicspath{{../}}
\usepackage{multirow}
\usepackage{float}
\usepackage[cmex10]{amsmath}
\usepackage {amssymb}
\usepackage{array}
\usepackage{mdwmath}
\usepackage{mdwtab}
\usepackage{eqparbox}
\usepackage{nomencl}
\usepackage{cite}
\usepackage{algorithm}
\usepackage{algorithmic}
\usepackage{makeidx}
\usepackage{ifthen}
\usepackage{subfigure}
\usepackage{caption}
\usepackage{pdflscape}
\usepackage{hyperref}
\usepackage{xcolor}

\newcommand{\beq}{\begin{equation}}
\newcommand{\eeq}{\end{equation}}
\newcommand{\beqn}{\begin{eqnarray}}
\newcommand{\eeqn}{\end{eqnarray}}
\newcommand{\beqno}{\begin{eqnarray*}}
\newcommand{\eeqno}{\end{eqnarray*}}
\newcommand{\bma}{\begin{displaymath}}
\newcommand{\ema}{\end{displaymath}}
\newcommand{\bnu}{\begin{enumerate}}
\newcommand{\enu}{\end{enumerate}}
\newcommand{\bce}{\begin{center}}
\newcommand{\ece}{\end{center}}
\newcommand{\btb}{\begin{tabular}}
\newcommand{\etb}{\end{tabular}}

\def\bx{{\mathbf{x}}}

\def\bq{{\mathbf{q}}}

\def\bP{{\mathbf{P}}}

\def\bc{{\mathbf{c}}}

\allowdisplaybreaks[4]

\begin{document}

\title{Optimal Workload Allocation for  Distributed Edge Clouds With Renewable Energy and Battery Storage}

\author{\IEEEauthorblockN{Duong~Thuy~Anh~Nguyen, Jiaming~Cheng, Ni~Trieu
and~Duong~Tung~Nguyen}  \vspace{-0.56cm}
\thanks{The first two authors have contributed equally to this work. 
}
}

 
 

\maketitle

\begin{abstract}
This paper studies an optimal workload allocation problem for a network of renewable energy-powered edge clouds that serve users located across various geographical areas. Specifically, each edge cloud is furnished with both an on-site renewable energy generation unit and a battery storage unit. Due to the discrepancy in electricity pricing and the diverse temporal-spatial characteristics of renewable energy generation, how to optimally allocate workload to different edge clouds to minimize the total operating cost while maximizing renewable energy utilization is a crucial and challenging problem. To this end, we introduce and formulate an optimization-based framework designed for Edge Service Providers (ESPs) with the overarching goal of simultaneously reducing energy costs and environmental impacts through the integration of renewable energy sources and battery storage systems, all while maintaining essential quality-of-service standards. Numerical results demonstrate the effectiveness of the proposed model and solution in maintaining service quality as well as reducing operational costs and emissions. Furthermore, the impacts of renewable energy generation and battery storage on optimal system operations are rigorously analyzed.
\end{abstract}

\printnomenclature


\vspace{-0.3cm}

\allowdisplaybreaks
\section{Introduction}
\label{Sec:Intro} 



Over the past decade, Cloud/Edge Service Providers (ESPs) have emerged as indispensable drivers of digital transformation.
Their pivotal role lies in facilitating the delivery of a wide spectrum of digital services, encompassing tasks such as data storage, processing, software applications, and beyond \cite{wiopt23}.
Each ESP oversees a portfolio of edge clouds (EC), also known as edge data centers. The ESPs typically manage extensive networks comprising variously sized and configured ECs, each housing a diverse array of computing resources.
The ECs are 
distributed across different geographical locations to ensure proximity to end-users, minimize latency, and optimize service delivery. 
These ECs serve as the fundamental building blocks of their cloud infrastructure, facilitating the provisioning of a diverse range of services to customers, ranging from virtual machines and storage to machine learning and content delivery. 

Energy efficiency is a critical consideration in cloud/edge computing. 
ECs usually consume a significant amount of energy to operate servers, networking equipment, cooling systems, and other infrastructure components. Most ECs are connected to the electrical grid and rely on utility power as their primary source of electricity. Due to their enormous energy consumption, ECs 
contribute to a large amount of greenhouse gas (GHG) emissions and increased carbon footprint. While the cost of electricity receives growing attention, the environmental impacts of power-intensive operations are often overlooked. Indeed, inexpensive electricity can sometimes come at the expense of environmental harm. According to the 2021 data from the U.S. Energy Information Administration \cite{EmissionData,ElectricityPriceData}, states like Wyoming, Utah, and North Dakota have some of the lowest electricity prices but significantly higher carbon footprints in their power sectors compared to the national average.

The exponential surge in data generation and computing demands has ushered in a relentless increase in the energy consumption of  ECs. This heightened energy demand poses both environmental and economic challenges. The commitment to utilizing more renewable energy sources, including solar, wind, and hydropower, has become one of the foremost strategic and operational goals for ECs, offering a solution to mitigate their carbon footprint, align with sustainability objectives, and decrease dependence on fossil fuels. A primary challenge associated with renewable energy sources is their intermittent nature. To address this, energy storage solutions like batteries are widely recognized as attractive options for promoting the sustainability and efficiency of EC operations. By storing surplus energy from renewable sources and/or low-cost grid electricity during off-peak hours to 
power ECs during peak periods or outages, they effectively tackle the variability and intermittency of renewable energy sources, facilitate their efficient utilization, and contribute to a resilient energy system. 

In regards to green EC design, \cite{LiangGreenDC2009} introduces an architecture for real-time monitoring, live virtual machine (VM) migration, and VM placement optimization to minimize power consumption. This approach aims to enhance server utilization and optimize power management in ECs, ultimately reducing their carbon footprint. Another line of research focuses on energy-cost-aware request routing among ECs by considering geographically dependent electricity costs \cite{Lei2010}. Recently, renewable energy resources 
have been integrated into ECs, thereby advancing sustainability 
\cite{Ghamkhari2013}. Reference \cite{Mohsenian2010} 
explores energy-information transmission trade-offs across various optimization problems, encompassing electricity costs, request routing, data center locations, proximity to renewable energy sources, server quantities, and the implications of carbon taxes. Reference \cite{Abbas2015} introduces the concepts of ``green workload" and ``green service rate" in contrast to ``brown workload" and ``brown service rate" which serves to distinctly address the separation of maximizing green energy utilization and minimizing brown energy costs. Reference \cite{hogade21} proposes a game theory-based resource management framework that incorporates renewable energy with the goal of minimizing cloud operating costs and queuing delays. In \cite{Zhang23}, a cooperative framework is considered where multiple electricity retailers work together to implement incentive-based demand response in distributed data centers, aiming to maximize profits. In\cite{jianwei_TSG19}, authors consider a workload allocation model for data center in electricity market.

\textbf{Contributions:}  Motivated by the compelling considerations outlined above, this paper proposes a holistic model for 
optimal EC operations. 
Given the diverse edge environments, 
characterized by varying electricity costs and carbon footprints, it becomes imperative for ESPs to formulate and implement efficient 
workload allocation strategies to ensure high quality-of-service (QoS), minimize 
costs, and enhance sustainability. Furthermore, the potential co-location of renewable energy generation and battery storage at EC facilities drives our investigation into tailored optimization models for co-optimizing EC provisioning and power procurement with these technologies. Specifically, we introduce an optimization framework for ESPs with the goal of simultaneously minimizing energy costs and environmental impact by integrating renewable energy sources and battery storage systems while maintaining essential QoS standards.

While our work shares some common ground with previous research, our core focus diverges notably. Unlike previous work that delves into decisions regarding EC placement, our primary objective centers on enhancing operational efficiency, guaranteeing the delivery of high-quality services while simultaneously curbing power consumption and reducing our environmental footprint. We employ a systematic and analytical framework, with a particular emphasis on environmental sustainability, wherein we meticulously account for emissions and carbon taxes. Our approach prominently features the integration of renewable energy sources and the implementation of battery storage as key components. Additionally, our model facilitates agile energy trading with the grid, offering a multifaceted solution that aligns with our overarching goals. Our primary focus is to analyze the impacts of integrating battery storage units, onsite renewable energy generation, and two-way energy trading with the grid on the optimal operation of networked ECs. Our numerical results illuminate the profound interdependencies among these factors and provide pragmatic insights for ESPs to reduce renewable energy curtailment. Previous work has often overlooked these aspects, focusing primarily on technical solutions to address intricate optimization models. 


\vspace{-0.3cm}

\section{System Model}
\label{Sec:SystemModel}

We consider an Edge Service Provider (ESP) that owns and operates multiple edge clouds (ECs) in different geographic locations. Each EC is equipped with servers, computing resources, and networking functions to provide low-latency, high-performance cloud services to users situated in various areas. The ESP operates with the overarching goal of delivering high-quality cloud services to its diverse user base while simultaneously optimizing several critical aspects of its operations, which encompass minimizing unmet demand, reducing energy consumption, and mitigating emissions. 

The ESP aggregates requests from users in different geographical regions and these requests are then directed to ECs for further processing. The ESP strives to ensure that users demand is consistently met, resulting in a seamless and reliable user experience. ECs typically consume a substantial amount of electrical energy to power servers, networking equipment, and other infrastructure components. Furthermore, due to their enormous energy consumption, ECs exert a substantial influence on the electric grid, contributing significantly to greenhouse gas (GHG) emissions and an overall increase in carbon footprint. Thus, energy efficiency is crucial to reduce operational costs and minimize the environmental impact of ECs' operations, with a particular focus on decreasing the carbon footprint. This involves deploying energy-efficient infrastructure components and strategically integrating renewable sources, like solar panels and wind turbines, along with battery storage systems to efficiently manage energy usage. By incorporating renewable energy sources and battery storage systems into EC facilities, surplus energy generated during off-peak hours can be stored and later utilized to power the EC during peak demand periods. Additionally, the ESP has the opportunity to sell excess power generated from its renewable sources back to the grid, generating extra revenue to offset energy costs and reduce both GHG emissions and dependence on fossil fuels.

\begin{figure}
    \centering
    \includegraphics[width=0.4\textwidth,height=0.16\textheight]{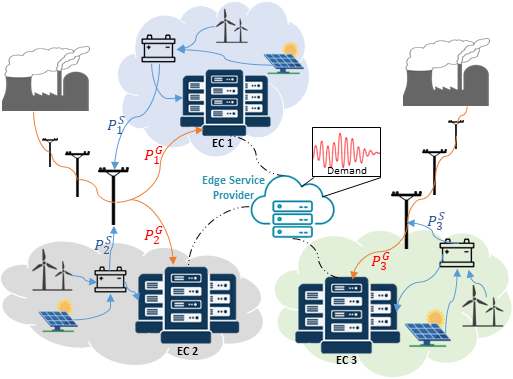}
    \caption{System model}
    \label{fig:system_model}
    \vspace{-0.5cm}
\end{figure}

To effectively govern its operations and deliver an exceptional user experience, the ESP adheres to operational constraints while upholding stringent QoS standards. These operational constraints encompass a wide range of factors, such as optimizing server activation based on varying demand at each EC, strategically distributing workloads to align with user needs, and maintaining strict thresholds for average delay and server utilization. This comprehensive approach ensures a consistently high level of QoS while concurrently reducing operational costs. Furthermore, the ESP faces additional constraints related to grid capacity limits and the availability of renewable energy sources. These constraints necessitate efficient energy management, including the fine-tuning of recharge and discharge rates for batteries to maintain balance and sustainability. By effectively balancing these operational and environmental constraints, the ESP aims to provide robust, sustainable, and environmentally responsible cloud services that meet the evolving demands and expectations of its diverse user base. The system model is depicted in Fig.~\ref{fig:system_model}.



    \vspace{-0.2cm}

\section{Problem Formulation}
\label{Sec:ProblemFormulation}
Consider a ESP with a total of $N$ ECs, denoted as $\mathcal{N}$, serving users situated in various areas. These geographical user areas are represented by Access Points (APs), collectively defined as the set $\mathcal{M}$, with $M$ individual APs identified using index $i$, while each EC is denoted by index $j$. In our model, we consider the time interval denoted as $\mathcal{T}$, where $\Delta T$ represents the duration of a single period within this interval.
\vspace{-0.3cm}

\subsection{Computing Capacity}
To enhance efficiency, optimize costs, and manage resources effectively, the ESP determines the number of servers to activate based on demand. Let $c_j^t$ represent the number of active computer servers at EC $j\in\mathcal{N}$ during time period $t \in \mathcal{T}$. It is imperative that the number of active servers does not surpass the total available servers at EC $j$, denoted as $C_j^{\max}$. This constraint is formally expressed as follows:
\begin{align}\label{eq-const-servers}
0 \le c_j^t \le C_j^{\max}.
\end{align}


\vspace{-0.9cm}
\subsection{Workload Allocation}
For each area $i \in \mathcal{M}$ and every time period $t \in \mathcal{T}$, we use $\lambda_i^t$ to represent the total expected resource demand of users in area $i$ at time $t$. Demand volumes can vary significantly throughout the day. 
We denote the number of requests from area $i$ allocated to EC $j$ during period $t$ as $x_{i,j}^t$, and $\alpha$ represents the average computing resource required to serve one request. Service requests from each area must either be served by ECs or considered unmet, denoted as $q_i^t$. Thus, we impose the following constraint:\vspace{-0.2cm}
\begin{align}\label{eq-const-WorkloadAllocation}
\alpha\left(\sum_{j}x_{i,j}^t+q_{i}^t\right) = \lambda_i^t.
\end{align}

Additionally, we assume a penalty cost of $\phi_i$ for each unmet request in area $i$. Therefore, the total unmet demand penalty for the ESP are calculated as:\vspace{-0.2cm}
\begin{align}\label{eq-cost-UnmetDemand}
C^{{\sf u}} = \sum_{i,t}\phi_i q_{i}^t.
\end{align}

To illustrate the impact of delay requirements on the system's performance, we denote $d_{i,j}$ as the propagation delay between each AP $i \in \mathcal{M}$ and EC $j \in \mathcal{N}$. To ensure Quality of Service (QoS), we introduce the binary indicator parameter $b_{i,j}$, which depends solely on the propagation delay $d_{i,j}$ and the maximum delay threshold $D^{\max}$ \cite{ResilientEC}. Specifically, the round-trip propagation delay must always be kept below $D^{\max}$. In other words,
\vspace{-0.1cm}
\begin{equation}
\label{delay}
   b_{i,j} = \left\{
\begin{array}{ll}
      1, & 2 d_{i,j} \leq D^{\max} \\
      0, & 2 d_{i,j} > D^{\max} \\      
\end{array}, ~~  \forall i,j. \right. 
\vspace{-0.05cm}
\end{equation}
\vspace{-1cm}



\subsection{Energy Consumption Model}
\subsubsection{Average Server Utilization}
Let $\rho$ represent the service rate signifying the maximum number of service requests that a single server can effectively handle or process within a time period. The average server utilization at EC $j$ during time period $t$, denoted as $\gamma_j^t$ and discussed in \cite{Mohsenian2010}, quantifies the proportion of the EC's capacity for that specific time period and is determined as follows:\vspace{-0.2cm}
\begin{align}\label{eq-AvgServerUtilization} 
    \gamma_j^t = \frac{\sum_{i\in\mathcal{M}}x_{i,j}^t}{\rho c_j^t}.
\end{align} 
To control queuing delay, we set a limit $\gamma^{\max} \in (0, 1]$ on the average server utilization at each EC as follows: \vspace{-0.2cm}
\begin{align}\label{eq-const-MaxAvgServerUtilization}
    \gamma_j^t \le \gamma^{\max}.
\end{align}
The choice of the $\gamma^{\max}$ parameter depends on the service request traffic pattern and the QoS requirements \cite{Mohsenian2010}. If $\gamma^{\max}$ is sufficiently small, the waiting time for a service request at an EC before server handling becomes negligible, with most of the overall latency in responding to service requests determined by the bounded propagation delay, $D^{\max}$.

\subsubsection{Power Consumption}
Let $P_j^{idle}$ denote the average power consumption of an individual server during idle periods, and let $P_j^{peak}$ denote the average power consumption when the server is actively processing service requests. Additionally, we introduce the term Power Usage Effectiveness (PUE)\footnote{PUE is a metric used to assess EC energy efficiency. As reported in the 2016 U.S. EC Energy Report, the average annualized PUE across various EC types typically falls within the range of $1.8$ to $1.9$.}, represented as $E_j^{usage}$.
The total power consumption (or power demand) $P_j^{U,t}$ at each EC location $j$ and for each period $t$ can be computed as follows \cite{Xiaobo2007,jianwei_TSG}:
\begin{align}\label{eq-PowerConsume}
    P_j^{U,t} &= c_j^t\left(P_j^{idle}+\left(E_j^{usage}-1\right)P_j^{peak}\right) \nonumber\\
    &+c_j^t\left(P_j^{idle}-P_j^{peak}\right)\gamma_j^t.
\end{align}


The ratio $P_j^{peak}/P_j^{idle}$ serves as a metric for assessing the power elasticity of servers. A higher value of this ratio indicates greater elasticity, resulting in reduced power consumption during periods of server inactivity. When $P_j^{peak}=P_j^{idle}$, the power consumption is  $P_j^{U,t} = c_j^tE_j^{usage}P_j^{peak}$. In this scenario, power consumption becomes solely dependent on the quantity of servers, without consideration for the number of routed requests or the operational period. 
\vspace{-0.3cm}
\subsection{Energy Model}
\subsubsection{Cost of Electricity}
ECs usually rely on the electrical grid as their primary power source. In North America, the electric grid operates on a regional basis. Most regions have regulated electricity markets where prices remain fixed throughout the day. In areas with deregulated markets, energy prices from the grid can fluctuate significantly throughout the day and across seasons, reflecting the dynamics of the wholesale electricity market. In our system model, we consider the general scenario where, at each EC $j$ and during each period $t$, the electricity price is represented as $e_j^t$. We introduce the variable $P_j^{G,t}$ to represent the amount of power to be purchased from the grid at a price $e_j^t$. Consequently, the corresponding electricity cost can be calculated as $e_j^t P_j^{G,t}.$


\subsubsection{Grid Capacity Limit}
The amount of power accessible for operating ECs at different geographic location is contingent on various factors such as the number of power plants in a region, their generation capacities especially those relying on renewable sources, and the existing residential, commercial, and industrial power demands. We take into consideration the ESP's maximum grid power draw at location $j$ during period $t$, which is constrained by $P_j^{G,\max,t}$, implying 
\begin{align}\label{eq-const-GridLimit}
    0 \le P_j^{G,t} \le P_j^{G,\max,t}.
\end{align}

\subsubsection{Renewable Energy}
To manage power costs effectively, reduce carbon emissions, and ensure uninterrupted operations, ESPs often employ a diverse range of energy sources and technologies. These strategies include integrating renewable energy sources like solar panels and wind turbines into their infrastructure. At time $t$, the EC operator has knowledge of the available renewable energy, denoted as $P_j^{R,t}$, at location $j$. 
If the renewable energy is sufficient to meet the current power demand, i.e., $P_j^{R,t} \ge P_j^{U,t}$, then no further procurement is necessary, meaning $P_j^{G,t} = 0$. Otherwise, the ESP must determine how much additional power ($P_j^{G,t}$) to purchase from the market.


\subsubsection{Battery Storage System}

Batteries with finite capacity can be incorporated to store excess energy during periods of low demand or when renewable sources produce surplus electricity, ensuring stable supply and lower costs. When no procurement is required, the surplus energy can be used for charging the battery. Conversely, if procurement is necessary, the ESP can decide whether to opportunistically discharge energy from the battery, in addition to the procured power, to meet the current demand.

Indeed, when surplus energy is stored in the batteries and not required to power ECs, the ESP also has the option to sell this excess energy back to the grid. Denoting the amount of electricity to be sold to the grid as $P_{j}^{S,t}$ at a sell-back price of $a_j^t$, we can calculate the electricity adjustment cost as follows:
\begin{align}\label{eq-cost-Electricity}
    C^{{\sf e}} = \sum_{j,t} \left(e_j^t P_j^{G,t} - a_j^t P_{j}^{S,t}\right).
\end{align}


\subsubsection{Power Balance Equation}
Let $P_j^{C,t}$ denote the charged battery energy and $P_j^{D,t}$ denote the amount of battery energy discharged. Constraints on maximum recharge and discharge rates for batteries are expressed as:
\begin{align}\label{eq-const-MaxRechargeDischarge}
    0 \le P_j^{C,t} \le P_j^{C,\max} \text{ and } 0 \le P_j^{D,t} \le P_j^{D,\max}.
\end{align} 

For simplicity, we are not considering power transmission losses. The energy balance equation is expressed as follows:
\begin{align}\label{eq-const-EnergyBalance}
\!P_j^{R,t}\!+\!P_j^{G,t}\!+\!P_j^{D,t}\!=\!P_j^{U,t}\!+\!P_j^{C,t}\!+\!P_j^{S,t}.
\end{align}
  

\subsubsection{Energy Dynamic}
We make the assumption of uniform efficiency for battery charging and discharging, represented by the parameter $\eta \in [0,1]$. Specifically, $\eta=0.8$ signifies that only $80\%$ of the energy is effectively utilized during both the charging and discharging processes. The dynamics governing the battery energy level at location $j$ and time $t$, denoted as $E_j^t$, can be expressed as follows:
\begin{align}\label{eq-const-EnergyLevel}
    E_j^{t+1}=E_j^t+\Delta T \left(\eta P_j^{C,t} - \frac{P_j^{D,t}}{\eta}\right),
\end{align}
where $\Delta T$ is the length of one time period.

\subsubsection{Battery Capacity Constraint}
Ensuring the desired EC availability is of utmost importance, especially when dealing with a finite battery capacity. We assume that the battery has a capacity of $E_j^{\max}$, and we choose to store a minimum energy of $E_j^{\min}$ to ensure EC availability, resulting in the following constraint: \vspace{-0.2cm}
\begin{align}\label{eq-const-BoundsEnergyLevel}
    E_j^{\min} \le E_j^t \le E_j^{\max}.
\end{align}

\subsection{Environmental Impact}
\subsubsection{Carbon Emission Factor}
The carbon emission factor for electricity varies significantly based on its source. In the case of electricity from the grid, the emission factor is location-dependent and closely tied to the energy composition of the region. Regions heavily reliant on fossil fuels, such as coal and natural gas, tend to exhibit higher carbon emission factors. Conversely, on-site renewable energy sources, such as solar panels and wind turbines, are renowned for their eco-friendliness. They typically boast significantly lower carbon emission factors or even approach zero emissions. 

In our model, we denote $ \theta_j$ as the emission factor associated with each unit of electricity purchased from the grid to power EC $j$. Consequently, the total carbon emissions linked to powering EC $j$ at time $t$ are calculated as follows:
\begin{align}\label{eq-Emission}
    EM_j^t = \theta_j P_j^{G,t}.
\end{align}

\subsubsection{Carbon Tax}
To address environmental concerns, several states in the United States and Canadian provinces have implemented carbon taxes \cite{carbontax}. Typically, these carbon taxes are imposed on power plants, which then transfer the cost of the carbon pricing to consumers through increased electricity prices. Consequently, environmental considerations are factored into our system model through the cost of electricity. However, given that carbon taxes are not yet widely adopted, we introduce them as distinct parameters in our research to better understand their impacts. To this end, we designate the carbon tax for each EC location $j$ as $\delta_j$, resulting in an additional cost of 
\begin{align}\label{eq-cost-CarbonCost}
    C^{{\sf c}} = \sum_{j,t} \delta_j EM_j^t = \sum_{j,t} \delta_j \theta_j P_j^{G,t} .
\end{align}

\vspace{-0.4cm}
\subsection{ESP Optimization Model}
The model for the ESP can be formulated as a Mixed-Integer Linear Programming (MILP) problem as follows:
\begin{align}
&\underset{\bP,\bx,\bq,\bc}{\text{minimize}} ~~~~~~ C^{{\sf u}} + C^{{\sf e}}+ C^{{\sf c}} \label{eq-Obj}\\
&\text{subject to}~~  
\eqref{eq-const-servers},\eqref{eq-const-WorkloadAllocation},\eqref{eq-const-MaxAvgServerUtilization},\eqref{eq-const-GridLimit},\eqref{eq-const-EnergyBalance}-\eqref{eq-const-BoundsEnergyLevel} 
\label{eq-ConstrSet}\\
&~~~~~~~~~~~~~~\!  P_j^{G,t},P_j^{C,t},P_j^{U,t},P_j^{D,t},P_j^{S,t}\ge 0,~ \forall j,t \nonumber\\
&~~~~~~~~~~~~~~\! 0 \leq x_{i,j}^t \leq b_{i,j} \lambda_i^{t} ,~ \forall i,j,t;~~ q_i^t\ge 0,c_j^t\in \mathbb{Z}^+, \forall j,t.~ \nonumber
\end{align}

The components of the objective function, as expressed in \eqref{eq-Obj}, encompass unmet demand penalties, electricity expenses, and costs associated with carbon dioxide emissions resulting from the operation of ECs. 
Complementing this objective function are a set of essential constraints, provided in \eqref{eq-ConstrSet}, that govern the ESP's operations. These constraints include considerations such as computing capacity, request routing, thresholds for average delay and server utilization, grid capacity limits, limits on power output from renewable energy sources, energy balance, recharge and discharge rates for batteries, energy dynamic equation and energy level limits to ensure the availability of ECs.

\section{Numerical Results}
\label{Sec:NumericalResults}
We consider an edge system comprising 8 edge clouds (ECs) ($N = 8$) and 10 APs ($M = 10$). The edge network topology is based on the cities and locations of randomly selected Equinix ECs \footnote{https://www.equinix.com/data-centers/americas-colocation, Access 2022.}. In the \textit{default setting}, we assume that all ECs are eligible to serve demand from every area, i.e., $b_{i,j}$ is set to $1$ for all $i$ and $j$. We will also perform sensitivity analysis on larger networks with more than $10$ areas. In our setup, we assume that all ECs are eligible to serve demand from every area. We consider $P^{idle}$ to be randomly generated from U$[0.45,0.55]$ kilowatt-hour (kWh), while $P^{peak}$ is randomly generated from U$[1.2,1.5]$ kWh. As reported in the 2016 U.S. EC energy report, we adopt $E^{usage}$ to fall within the interval of U$[1.8,1.9]$. For the prices $e_j$ associated with selected locations, which range from $[0.1,0.35]$ per kWh, we obtain data from \cite{ElectricityPriceData}. Due to the absence of carbon tax data in certain U.S. states, we generate $\delta_j$ randomly from the range $[0.6,0.7]$ \cite{carbontax}.

The ``sell-back" price is expected to be less than the procurement cost, as denoted by $a_j = \zeta e_j, \forall j$, with $\zeta$ equal to $0.8$. By utilizing the trace data from GWA-DAS \footnote{http://gwa.ewi.tudelft.nl/datasets/gwa-t-1-das2/report}, we randomly generate the expected demand $\lambda_i^t$, ranging from $10$ to $30$ per hour. In our problem, we assume that each EC is directly connected to the grid. Thus, the grid capacity at location $j$ is randomly generated from U$[1,1.5]$ megawatts (MW). Regarding the emission factor for each EC $j$, we will focus solely on CO2 emissions. According to \cite{carbontax}, we assume that the carbon tax for the selected location $j$ follows U$[20,50]$ per ton of CO2 emission. The carbon emission factor $\theta_j$ is randomly generated from U$[0.1,0.8]$, based on data from \cite{EmissionData}.

In the \textit{default setting},  we consider that $P_{j,t}^{R}$ is randomly generated from U$[80,100]$ kW. The maximum battery capacity $E_j^{\max}$ is generated from U$[90,100]$ kW, while $E_j^{\min}$ follows U$[30,50]$ kW. The charging capacity ($P_{j}^{C,\max}$) and discharging capacity ($P_{j}^{D,\max}$) are generated from U$[70,80]$ kW and U$[70, 80]$ kW, respectively. Additionally, we consider the values $\gamma^{\max} = 0.9$, $\alpha = 0.5$, $\rho = 0.8$, $\Delta T = 1$, and $T = 12$. We will also vary these parameters during sensitivity analysis. All the experiments are conducted in MATLAB using CVX\footnote{http://cvxr.com/cvx/}  and Gurobi\footnote{https://www.gurobi.com/} on a desktop with an Intel Core i7-11700KF CPU and 32GB of RAM. 

\vspace{-0.3cm}
\subsection{Sensitivity analysis}
 This section presents a sensitivity analysis to assess the influence of key system parameters on the optimal solution. These parameters include renewable energy ($P^{R}{j}$), electricity price $e_j^t$, and the ``sell-back'' ratio ($\zeta$). To evaluate the impact of renewable energy on system performance, we introduce a scaling factor $\Psi$ for $P^{R}{j}$, where $\Psi = 1$ indicates the default value. Specifically, the value of $P^{R}{j}$ generated in the default setting is multiplied by $\Psi$ to either scale up or down the renewable energy. Similarly, we use $\xi{E^{\max}}, \xi_e,$ and $\xi_{D^{\max}}$ as scaling factors for the maximum battery size at EC ($E_j^{\max}$), electricity price $e_j$, and average utilization $D^{\max}$. In this sensitivity analysis, we focus exclusively on the proposed model.

\begin{figure}[h!]
		 \subfigure[varying $\Psi$ and $e$]{
	     \includegraphics[width=0.242\textwidth,height=0.09\textheight]{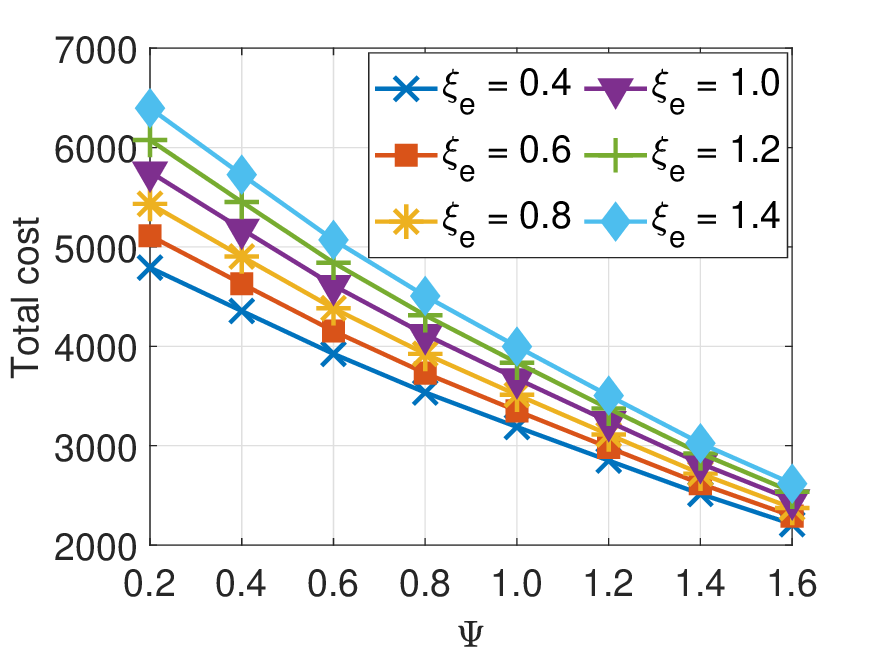}
	     \label{fig:Psi_e}
	}   \hspace*{-2.1em} 
	    \subfigure[varying $e$ and $E^{\max}$]{
	     \includegraphics[width=0.242\textwidth,height=0.09\textheight]{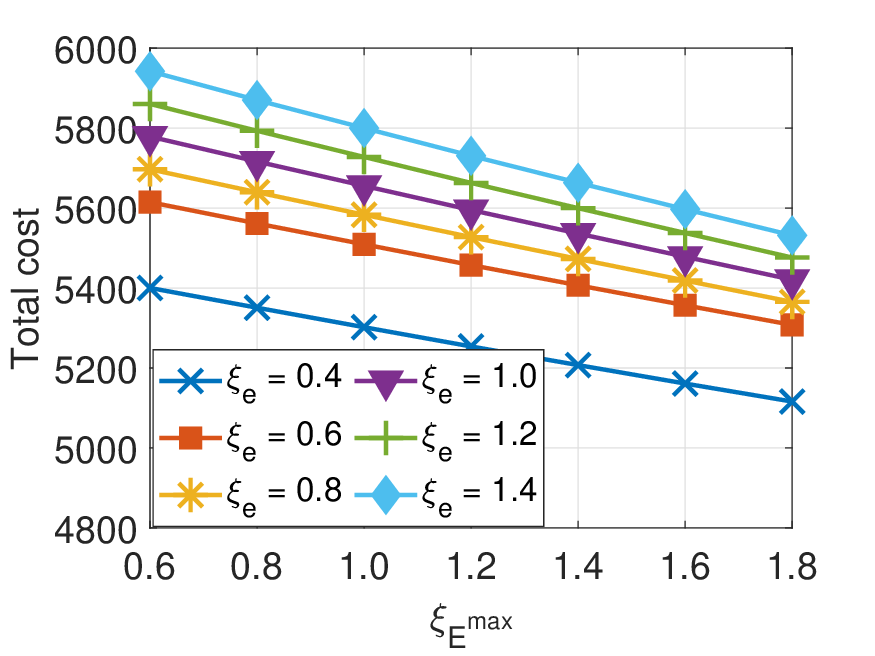}
	     \label{fig:Psi_E}
	} 	    \vspace{-0.2cm}
	    \caption{The benefits of renewable energy}
	    \vspace{-0.4cm}
\end{figure}

\noindent \textbf{\textit{1) Benefits of battery storage:}}
As depicted in Figure \ref{fig:Psi_e}, the total cost decreases as $\Psi$ increases, signifying an increase in the available renewable energy at each EC. This allows the operator to have greater flexibility in supplying power from renewable sources, resulting in a reduced reliance on grid energy. Furthermore, it is evident that a higher electricity price ($e$) can motivate the operator to maximize the utilization of renewable energy resources, reducing the need for energy procurement from the grid. Similarly, we also examine the impact of battery capacity size on EC operations. Recall that $\xi_{E^{\max}}$ is a scaling factor for battery size. When $\xi_{E^{\max}}$ is set to a higher value, each EC can potentially have a greater capacity to absorb surplus energy. This advantage becomes particularly significant when electricity prices ($e$) are elevated.

\begin{figure}[h!]
		 \subfigure[Revenue: varying $\Psi$ and $\zeta$]{
	     \includegraphics[width=0.242\textwidth,height=0.09\textheight]{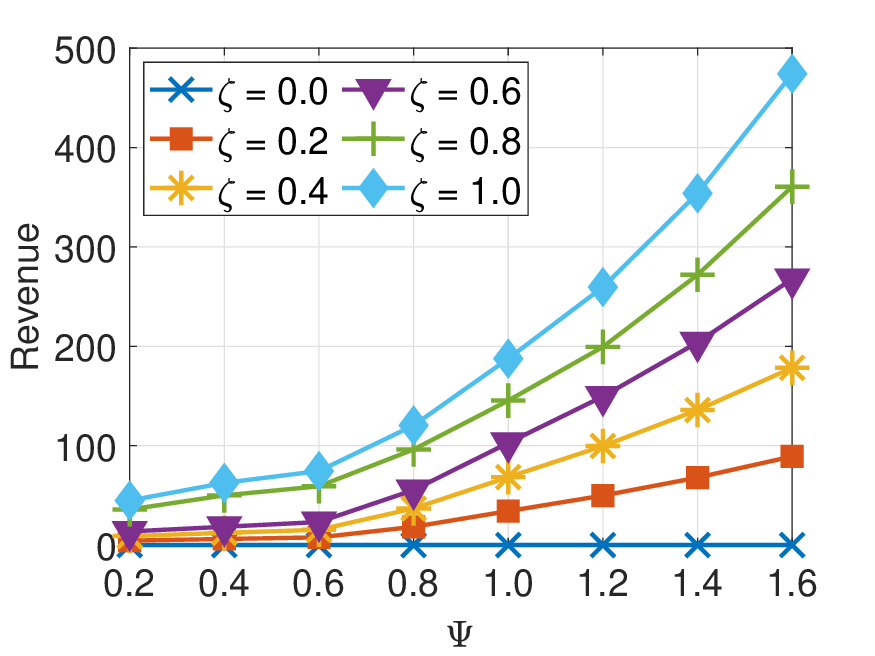}
	     \label{fig:zeta_Psi_rev}
	}   \hspace*{-2.1em} 
	    \subfigure[Total cost: varying $\Psi$ and $\zeta$]{
	     \includegraphics[width=0.242\textwidth,height=0.09\textheight]{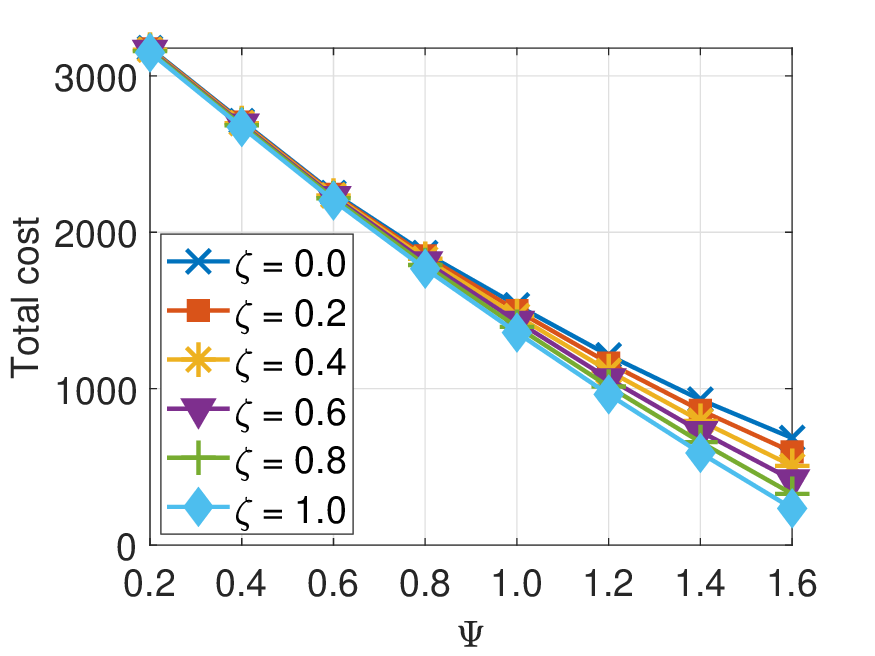}
	     \label{fig:zeta_Psi_cost}
	} 
 	 \subfigure[Revenue: varying $E^{\max}$ and $\zeta$]{
	     \includegraphics[width=0.242\textwidth,height=0.09\textheight]{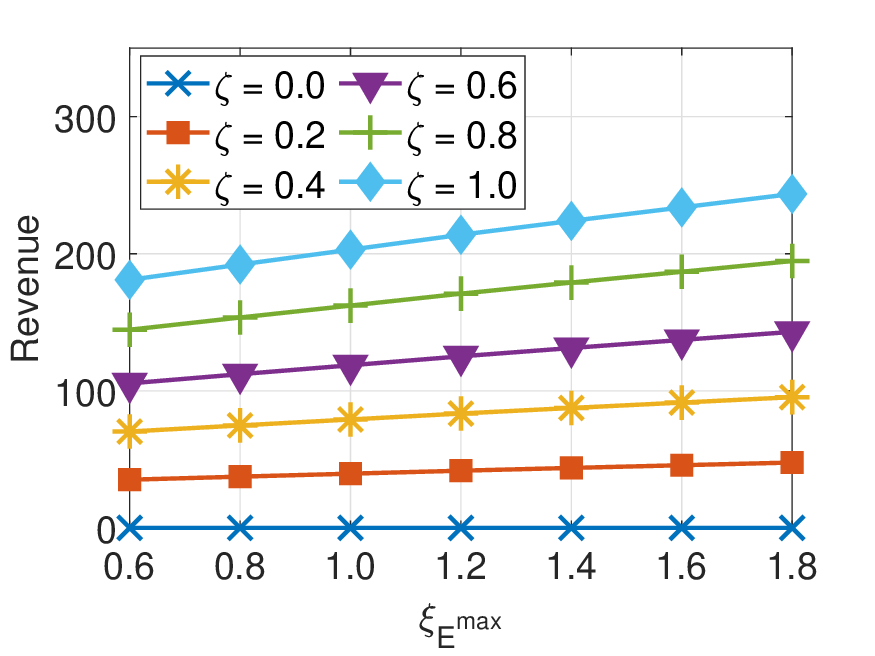}
	     \label{fig:zeta_E_rev}
	}   \hspace*{-2.1em} 
	    \subfigure[Total cost: varying $E^{\max}$ and $\zeta$]{
	     \includegraphics[width=0.242\textwidth,height=0.09\textheight]{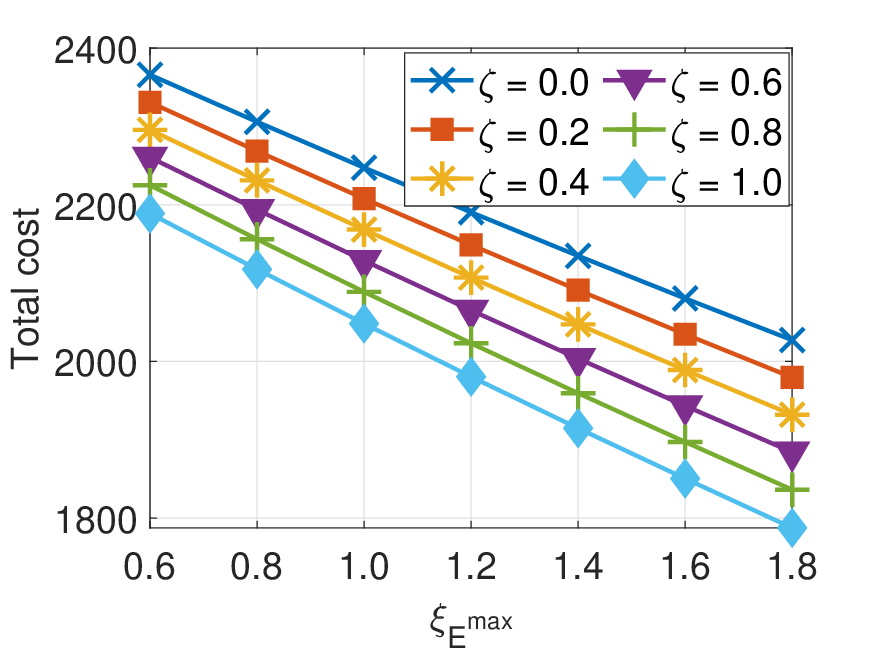}
	     \label{fig:zeta_E_cost}
	} 	    \vspace{-0.2cm}
	    \caption{The benefits of ``sell-back" option}
	    \vspace{-0.4cm}
\end{figure}

\noindent \textbf{\textit{2) Benefits of ``sell-back" option:}} Fig.\ref{fig:zeta_Psi_rev} - \ref{fig:zeta_E_cost} shows how ``sell-back ratio" influences the system performance. Recall that $\zeta$ is defined as ``sell-back" ratio between the electricity price and ``sell-back" price (i.e., $e_j = \zeta a_j, \forall j$). Once the power demand of each EC has been met, any surplus renewable energy tends to be given higher priority for selling back to the grid, especially when $\zeta$ is set to higher values. As illustrated in Fig.\ref{fig:zeta_E_rev} and \ref{fig:zeta_E_cost}, larger values of $E$ signify a greater battery capacity at each EC, allowing for more energy storage. In such cases, the operator may opt to store excess energy and subsequently sell it back to the grid, particularly when the ``sell-back" price is favorable.
\vspace{-0.4cm}

\begin{figure}[h!]
		 \subfigure[varying $\gamma^{\max}$ and $D^{\max}$]{
	     \includegraphics[width=0.242\textwidth,height=0.09\textheight]{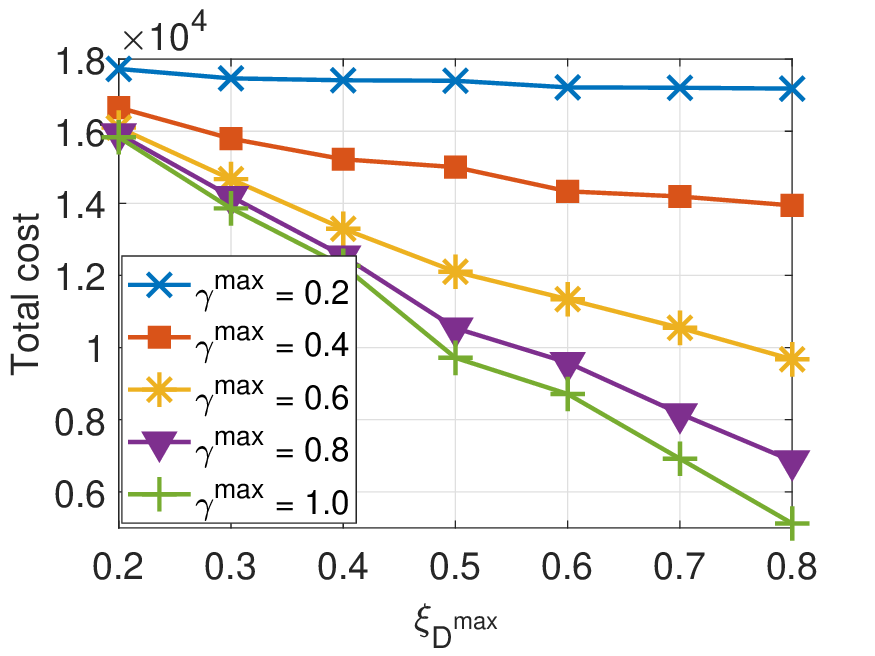}
	     \label{fig:gamma_b}
	}   \hspace*{-2.1em} 
	    \subfigure[varying $M$ and $N$]{
	     \includegraphics[width=0.242\textwidth,height=0.09\textheight]{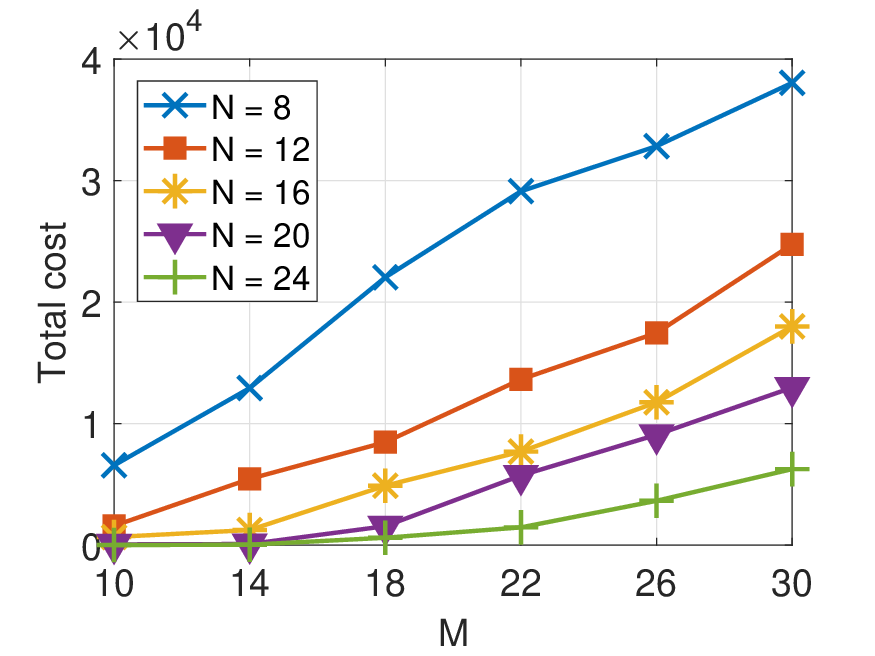}
	     \label{fig:EC_AP}
	} 	    \vspace{-0.2cm}
	    \caption{The impact of other system parameters}
     	\vspace{-0.6cm}
\end{figure}

\noindent \textbf{\textit{3) The impact of other system parameters:}}
Notice that $\gamma^{\max}$ and $D^{\max}$ directly influence QoS. Specifically, $\gamma^{\max}$ imposes a constraint on the average server utilization at each EC, while $D^{\max}$ sets a limit on the propagation delay between APs and ECs. As shown in Figures \ref{fig:gamma_b}, when $\omega_\gamma$ is decreased, it signifies a stricter restriction on the average utilization at each EC, which can result in a higher chance of unmet demand due to these more stringent requirements. As described in (\ref{delay}), an EC $j$ can only serve user requests from area $i$ when the propagation delay between them is within the threshold $D^{\max}$. In other words, a decrease in $D^{\max}$, indicating a more stringent delay requirement, leads to an increase in the number of $b_{i,j}$ values that become zero, indicating a reduction in the number of eligible ECs to serve user demand. Consequently, the total cost of the system increases as $D^{\max}$ decreases. Furthermore, Figure \ref{fig:EC_AP} shows the impact of network size on the optimal solution. As expected, with a fixed number of ECs, the total cost rises as the number of APs increases.

\vspace{-0.4cm}

\subsection{Performance Comparisons}
In this section, we aim to compare the performance of the proposed model with the following benchmarks:
\begin{itemize}
    \item \textbf{M1}: This model lacks the capability for the ``sell-back-to-grid" option, and ECs do not come equipped with batteries.
    \item \textbf{M2}: This model exclusively focuses on using batteries to store energy and does not allow the operator to sell surplus energy back to the grid.
    \item \textbf{M3}: This model is designed to enable the selling of excess energy back to the grid, while ECs do not feature battery installations.
\end{itemize}

The evaluation and comparison of the four schemes are based on their total cost with four different settings. To simplify, we refer to our proposed model as ``\textbf{M0}". These four models can be straightforwardly categorized into those that incorporate or omit considerations for sell-back and battery storage options.
 
 \vspace{-0.4cm}
\begin{figure}[h!]
		 \subfigure[Total cost: varying $\zeta$]{
	     \includegraphics[width=0.240\textwidth,height=0.09\textheight]{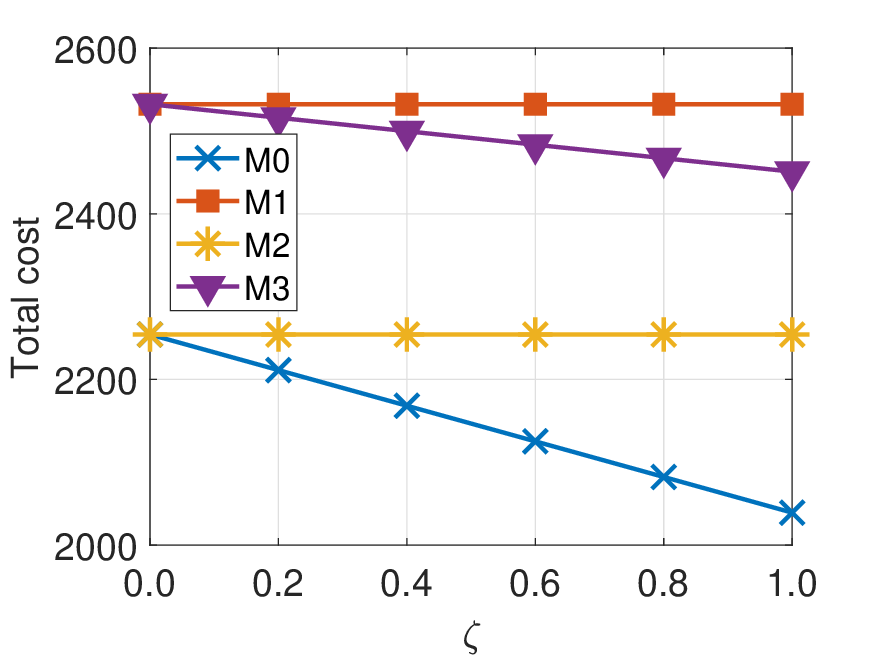}
	     \label{fig:zeta_cost_compare}
	}   \hspace*{-2.1em} 
 \vspace{-0.2cm}
	    \subfigure[Revenue: varying $\zeta$]{
	     \includegraphics[width=0.240\textwidth,height=0.09\textheight]{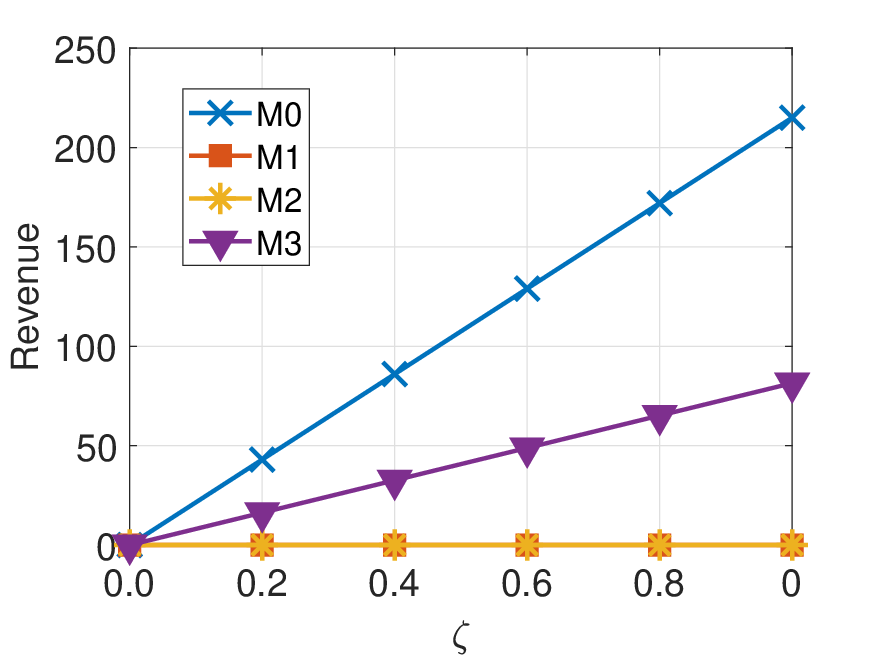}
	     \label{fig:zeta_rev_compare}
	} 
 	\subfigure[Total cost: varying $\xi_e$]{
	     \includegraphics[width=0.240\textwidth,height=0.09\textheight]{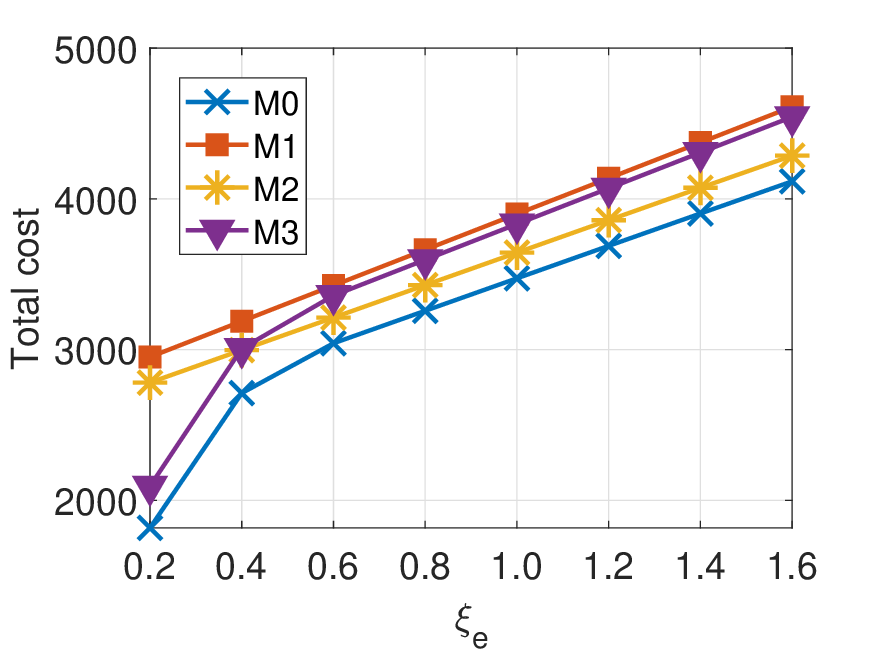}
	     \label{fig:e_compare}
	} \hspace*{-2.1em}
 	\subfigure[Total cost: varying $\xi_{E^{\max}}$]{
	     \includegraphics[width=0.240\textwidth,height=0.09\textheight]{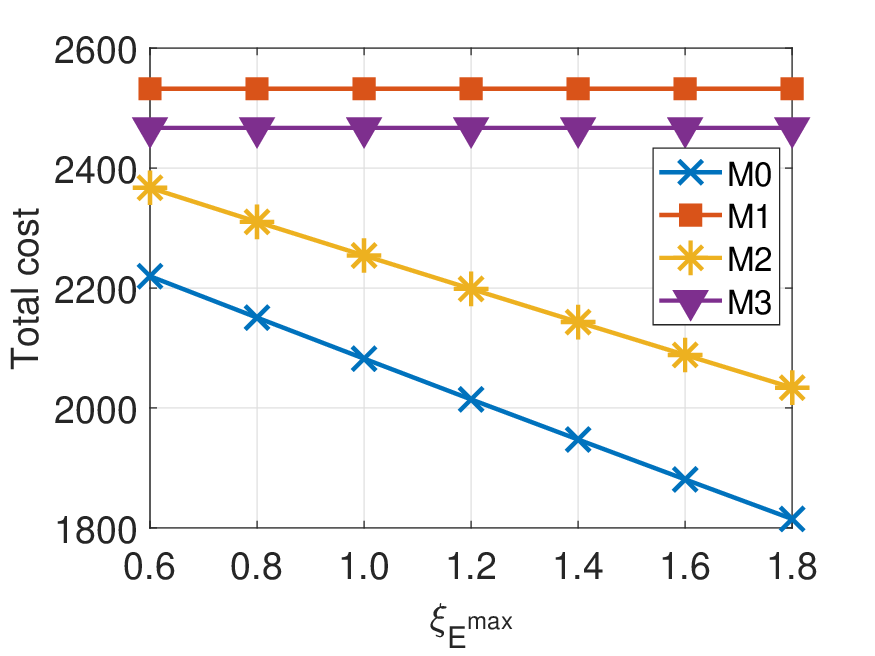}
	     \label{fig:E_max_compare}
	} 	    \vspace{-0.2cm}
	    \caption{Model comparisons}
     \vspace{-0.2cm}
\end{figure}

As illustrated in Fig. \ref{fig:zeta_cost_compare}--\ref{fig:E_max_compare}, our proposed model significantly outperforms the other schemes. \textbf{M1} achieves the worst performance since it lacks both the ``sell-back-to-grid" option and battery storage. Fig. \ref{fig:zeta_cost_compare} shows that \textbf{M0} and \textbf{M3} do not vary regardless of changes in $\zeta$ because these two models lack considerations for ``sell-back-to-grid". Thus, when $\zeta = 0$, indicating a sell-back price of $0$, there is no distinction in performance between \textbf{M0} and \textbf{M2} or \textbf{M1} and \textbf{M3}. The difference between these two pairs demonstrates the advantages of battery storage, as energy can be discharged as one of the available renewable resources. When $\zeta$ increases, the advantages of models that consider ``sell-back-to-grid'' options (such as \textbf{M0} and \textbf{M3}) become increasingly pronounced due to the elevated sell-back prices, which can be verified in Fig. \ref{fig:zeta_rev_compare}.

Furthermore, as depicted in Fig. \ref{fig:e_compare}, all four models exhibit an increase as electricity prices ($e$) rise. Notably, the advantages of the sell-back option can be emphasized, especially when electricity prices ($e$) are low. In such scenarios, selecting energy supply from the grid is not that expensive, offering greater operational flexibility. However, as the electricity price ($e$) escalates, a model with solely a sell-back option becomes insufficient, as these resources would be wasteful without storing surplus energy in a battery. Therefore, the cost of \textbf{M2} gradually surpasses that of \textbf{M3}. The reduction in the maximum battery size ($E^{\max}$) signifies decreased capacity for storing renewable energy, potentially resulting in limited availability of renewable energy and, consequently, the possibility of renewable energy curtailment. Figure \ref{fig:E_max_compare} illustrates that costs decrease across models that take battery storage into account, as the operator can prioritize the utilization of renewable energy sources by discharging energy from the batteries to meet power demand. These advantages become more prominent as the maximum battery size ($E^{\max}$) increases. In summary, the risk of renewable energy curtailment is heightened when EC systems lack both battery storage and the sell-back option.

\vspace{-0.3cm}

\section{Conclusion}
\label{conc}
This paper has unveiled valuable insights and potential benefits of integrating renewable energy sources, battery systems and the practice of selling back surplus energy in edge service operations. Notably, when renewable energy sources are considered, our approach underscores the allocation of workload not only to ECs with low electricity prices but also to those with high renewable generation capacity. Moreover, the introduction of batteries adds another dimension to the workload allocation problem as the battery capacity can efficiently absorb surplus energy, acting as an essential tool for energy arbitrage. Our future work will investigate flexible workload scheduling, enabling task deferral to accommodate varying energy availability, aligning with demand response principles. Another exciting direction is to consider and integrate various system uncertainties into the operational model of the networked ECs.

\vspace{-0.4cm}

\bibliographystyle{IEEEtran}
\bibliography{ref.bib}

\end{document}